# On the periodic continued radicals of 2 and generalization for Vieta's product

**Jayantha Senadheera**   jayantha.senadheera@gmail.com

**Abstract**

In this paper we study periodic continued radicals of 2. We show that any periodic continued radicals of *2* converge to $2\sin q\pi$, for some rational number *q* depends on the continued radical. Furthermore we show that if $r_n$ is a periodic nested radicals of *2*, which has *n* nested roots, then the limit points of the sequence $\langle 2^n(2\sin q\pi - r_n)\rangle$ have the form $\alpha\pi$, where $\alpha$ is an algebraic number. This result give a set of subsequences converges to $\alpha\pi$, for each $\alpha$. Also we show that limit of these sub sequences can be represented as "Vieta like" nested radical products. Hence this result generalizes the Vieta's product for π. Several interesting examples are illustrated.

**Introduction**

One objective of this paper is to establish basic theorems to understand the behavior and structure of continued radicals (infinite nested radicals) of 2. The other objective is to generalize well known Vieta's product for $\pi$. In mathematics literature we can find several attempts related to the first objective. Among these, problems 176, 177, 178 of Polya and Zego's excellent book [1] is the prominent one. In the past history, time to time several authors concerned about this topic [2],[3]. Edouard Lucas used nested radicals of 2 to develop his methods for primality testing [4]. Also the recent papers[5],[6],[7] and [8] are on these lines. In this paper, section 1 concentrates on the first objective. Most of the results in section 1 can be found in above references. But for the sake of completeness and for some generalizations, we reprove those results. At the end of section 1, we give a simple criterion to classify all periodic continued radicals of 2.

Basically section 2 concentrates on the second objective. In some sense it is the natural extension of section 1. Surprisingly this natural extension beautifully connects with generalization of Vieta's product. This generalization restricted to the periodic continued radicals of 2. For the different kind of generalizations, specially related with elliptic functions, refer the papers [9], [10].

In 1593 French mathematician Francois Vieta gave the beautiful infinite product formula for $\pi$ using only 2's and square roots.

$$\frac{2}{\pi} = \frac{\sqrt{2}}{2} \cdot \frac{\sqrt{2+\sqrt{2}}}{2} \cdot \frac{\sqrt{2+\sqrt{2+\sqrt{2}}}}{2} \cdots \qquad (1)$$

It is said that this was the first exact analytic expression given for $\pi$, as well as the first recorded use of an infinite product in Mathematics.

In his paper [9], Aaron Levin writes "Given the simplicity, elegance, and age of Vieta's product it is surprising that there seems to have been few attempts at finding similar formulas." In that paper he gave the following Vieta like product for $\frac{2\sqrt{3}}{\pi}$.

$$\frac{2\sqrt{3}}{\pi} = \frac{1+\sqrt{2}}{2} \cdot \frac{1+\sqrt{2-\sqrt{2}}}{2} \cdot \frac{1+\sqrt{2-\sqrt{2-\sqrt{2}}}}{2} \cdot \ldots \qquad (2)$$

Certainly it is a surprise that it took more than four hundred years to come with this second Vieta like product formula from the appearance of first one. By seeming these two formulas one can suspect that whether there is a more general setting behind these products. In this paper we give an elementary generalization for Vieta like products. In particular in this generalization we give perhaps the next simplest expression of Vieta like product family.

$$\frac{(\sqrt{5}-1)5\sqrt{2}}{\pi\sqrt{5-\sqrt{5}}} = \frac{(\sqrt{5}+1)\frac{1}{2}+\sqrt{2}}{2} \left( \frac{(\sqrt{5}-1)\frac{1}{2}+\sqrt{2-\sqrt{2}}}{2} \cdot \frac{(\sqrt{5}+1)\frac{1}{2}+\sqrt{2+\sqrt{2-\sqrt{2}}}}{2} \right) \left( \frac{(\sqrt{5}-1)\frac{1}{2}+\sqrt{2-\sqrt{2+\sqrt{2-\sqrt{2}}}}}{2} \cdot \frac{(\sqrt{5}+1)\frac{1}{2}+\sqrt{2+\sqrt{2-\sqrt{2+\sqrt{2-\sqrt{2}}}}}}{2} \right) \ldots \qquad (3)$$

## Section 1

Let $(\varepsilon_n)$ be a sequence such that $\varepsilon_n \in \{1, -1\}$ for all $n \in \mathbb{N}$. i.e $(\varepsilon_n) \in \{-1, 1\}^{\mathbb{N}}$. Let the sequence $(r_n)$ be defined as $r_n = \sqrt{2 + \varepsilon_1\sqrt{2 + \ldots + \varepsilon_n\sqrt{2}}}$ with $r_0 = \sqrt{2}$. Let $(c_n)$ be any real sequence such that $c_n \in [-2, 2]$ for all $n \in \mathbb{N} \cup \{0\}$. Define the sequence $(r_n(c_n))$ as follows.

$r_n(c_n) = \sqrt{2 + \varepsilon_1\sqrt{2 + \ldots + \varepsilon_n\sqrt{2 + c_n}}}$ with $r_0(c_0) = \sqrt{2 + c_0}$. Then observe that the sequence $r_n(c_n)$ is real for all $n \in \mathbb{N}$.

**Theorem 1**

The sequence $(r_n(c_n))$ convergent and $\lim_{n \to \infty} r_n(c_n) = 2\sin\left(1 + \sum_{i=1}^{\infty}\frac{\varepsilon_1\varepsilon_2\ldots\varepsilon_i}{2^i}\right)\frac{\pi}{4}$.

First we prove the following lemma.

**Lemma 1.1**

(i) $r_n(0) = r_n = 2\sin\left(1 + \sum_{i=1}^{n}\frac{\varepsilon_1\varepsilon_2\ldots\varepsilon_i}{2^i}\right)\frac{\pi}{4}$ and $r_n(2) = 2\sin\left(1 + \sum_{i=1}^{n}\frac{\varepsilon_1\varepsilon_2\ldots\varepsilon_i}{2^i} + \frac{\varepsilon_1\varepsilon_2\ldots\varepsilon_n}{2^n}\right)\frac{\pi}{4}$.

(ii) For any $n \in \mathbb{N}$, $r_{n-1}(0) \leq r_n(c_n) \leq r_n(2)$ or $r_{n-1}(0) \geq r_n(c_n) \geq r_n(2)$.

Proof of Lemma1.1

To show $r_n = 2\sin\left(1 + \sum_{i=1}^{n} \frac{\varepsilon_1 \varepsilon_2 ... \varepsilon_i}{2^i}\right)\frac{\pi}{4}$, observe that, when $n = 1$ the right hand side of the identity is $2\sin\left(1 + \frac{\varepsilon_1}{2}\right)\frac{\pi}{4}$.

But $2\sin\left(1 + \frac{\varepsilon_1}{2}\right)\frac{\pi}{4} = 2\sin\left(\frac{\pi}{4} + \frac{\varepsilon_1 \pi}{8}\right) = \sqrt{2 - 2\cos\left(\frac{\pi}{2} + \frac{\varepsilon_1 \pi}{4}\right)} = \sqrt{2 + 2\sin\frac{\varepsilon_1 \pi}{4}} = \sqrt{2 + \varepsilon_1 \sqrt{2}} = r_1$.

Hence when $n = 1$ the result is true. (Observe that we choose positive square root since

$0 \leq \frac{\pi}{4} + \frac{\varepsilon_1 \pi}{8} \leq \frac{\pi}{2}$. Therefore *sin* value is positive)

Suppose for $n \in \mathbb{N}$ the result is true. Observe that

$$2\sin\left(1 + \sum_{i=1}^{n+1} \frac{\varepsilon_1 \varepsilon_2 ... \varepsilon_i}{2^i}\right)\frac{\pi}{4} = \sqrt{2 - 2\cos\left(\frac{\pi}{2} + \sum_{i=1}^{n+1} \frac{\varepsilon_1 \varepsilon_2 ... \varepsilon_i}{2^{i-1}}\frac{\pi}{4}\right)}^{*}$$

$$= \sqrt{2 + 2\sin\left(\sum_{i=1}^{n+1} \frac{\varepsilon_1 \varepsilon_2 ... \varepsilon_i}{2^{i-1}}\right)\frac{\pi}{4}} = \sqrt{2 + 2\sin\left(\varepsilon_1 + \varepsilon_1 \sum_{i=1}^{n} \frac{\varepsilon_{1+1}\varepsilon_{1+2}...\varepsilon_{1+i}}{2^i}\right)\frac{\pi}{4}}$$

$$= \sqrt{2 + \varepsilon_1 2\sin\left(1 + \sum_{i=1}^{n} \frac{\varepsilon_{1+1}\varepsilon_{1+2}...\varepsilon_{1+i}}{2^i}\right)\frac{\pi}{4}} = \sqrt{2 + \varepsilon_1 \sqrt{2 + ... + \varepsilon_{n+1}\sqrt{2}}}.$$

Hence by induction result is true.

* Again we choose positive square root since $0 \leq \frac{\pi}{4} + \sum_{i=1}^{n+1} \frac{\varepsilon_1 \varepsilon_2 ... \varepsilon_i}{2^i}\frac{\pi}{4} \leq \frac{\pi}{2}$.

To show $r_n(2) = 2\sin\left(1 + \sum_{i=1}^{n} \frac{\varepsilon_1 \varepsilon_2 ... \varepsilon_i}{2^i} + \frac{\varepsilon_1 \varepsilon_2 ... \varepsilon_n}{2^n}\right)\frac{\pi}{4}$, observe that, when $n = 1$ the right hand side of the identity is $2\sin(1 + \varepsilon_1)\frac{\pi}{4} = \sqrt{2 - 2\cos\left(\frac{\pi}{2} + \frac{\varepsilon_1 \pi}{2}\right)} = \sqrt{2 + 2\sin\left(\frac{\varepsilon_1 \pi}{2}\right)} = \sqrt{2 + \varepsilon_1 \sqrt{2 + 2}} = r_1(2)$.

Hence when $n = 1$ the result is true. Suppose for $n \in \mathbb{N}$ the result is true. Observe that

$$2\sin\left(1 + \sum_{i=1}^{n+1} \frac{\varepsilon_1 \varepsilon_2 ... \varepsilon_i}{2^i} + \frac{\varepsilon_1 \varepsilon_2 ... \varepsilon_{n+1}}{2^{n+1}}\right)\frac{\pi}{4} = \sqrt{2 - 2\cos\left(\frac{\pi}{2} + \sum_{i=1}^{n+1} \frac{\varepsilon_1 \varepsilon_2 ... \varepsilon_i}{2^{i-1}}\frac{\pi}{4} + \frac{\varepsilon_1 \varepsilon_2 ... \varepsilon_{n+1}}{2^n}\frac{\pi}{4}\right)}$$

$$= \sqrt{2+2\sin\left(\sum_{i=1}^{n+1}\frac{\varepsilon_1\varepsilon_2...\varepsilon_i}{2^{i-1}}+\frac{\varepsilon_1\varepsilon_2...\varepsilon_{n+1}}{2^n}\right)\frac{\pi}{4}} = \sqrt{2+2\sin\left(\varepsilon_1+\varepsilon_1\sum_{i=1}^{n}\frac{\varepsilon_{1+1}\varepsilon_{1+2}...\varepsilon_{1+i}}{2^i}+\frac{\varepsilon_{1+1}\varepsilon_{1+2}...\varepsilon_{1+n}}{2^n}\right)\frac{\pi}{4}}$$

$$= \sqrt{2+\varepsilon_1 2\sin\left(1+\sum_{i=1}^{n}\frac{\varepsilon_{1+1}\varepsilon_{1+2}...\varepsilon_{1+i}}{2^i}+\frac{\varepsilon_{1+1}\varepsilon_{1+2}...\varepsilon_{1+n}}{2^n}\right)\frac{\pi}{4}} = \sqrt{2+\varepsilon_1\sqrt{2+...+\varepsilon_n\sqrt{2+\varepsilon_{n+1}\sqrt{2+2}}}} = r_{n+1}(2)$$

Hence the result.

Proof of (ii):-

Observe that since $-2 \leq c_1 \leq 2$, $0 \leq \sqrt{2+c_1} \leq 2 \Rightarrow 0 \leq \varepsilon_1\sqrt{2+c_1} \leq 2\varepsilon_1$ or $0 \geq \varepsilon_1\sqrt{2+c_1} \geq 2\varepsilon_1$

$\Rightarrow \sqrt{2} \leq \sqrt{2+\varepsilon_1\sqrt{2+c_1}} \leq \sqrt{2+2\varepsilon_1}$ or $\Rightarrow \sqrt{2} \geq \sqrt{2+\varepsilon_1\sqrt{2+c_1}} \geq \sqrt{2+2\varepsilon_1}$

$\Rightarrow r_0(0) \leq r_1(c_1) \leq r_1(2)$ or $r_0(0) \geq r_1(c_1) \geq r_1(2)$. Therefore when $n=1$ the result is true.

Now suppose for $n \in \mathbb{N}$ the result is true. i.e $r_{n-1}(0) \leq r_n(c_n) \leq r_n(2)$ or $r_{n-1}(0) \geq r_n(c_n) \geq r_n(2)$.

Observe that $r_{n+1}(c_{n+1})$ can be written as $r_{n+1}(c_{n+1}) = \sqrt{2+r_n'(c_{n+1})}$,

where $r_n'(c_{n+1}) = \sqrt{2+\varepsilon_{1+1}\sqrt{2+...+\varepsilon_{1+n}\sqrt{2+c_{n+1}}}}$.

But by induction hypothesis $r_{n-1}'(0) \leq r_n'(c_{n+1}) \leq r_n'(2)$ or $r_{n-1}'(0) \geq r_n'(c_{n+1}) \geq r_n'(2)$.

$\therefore \sqrt{2+\varepsilon_1 r_{n-1}'(0)} \leq \sqrt{2+\varepsilon_1 r_n'(c_{n+1})} \leq \sqrt{2+\varepsilon_1 r_n'(2)}$ or $\sqrt{2+\varepsilon_1 r_{n-1}'(0)} \geq \sqrt{2+\varepsilon_1 r_n'(c_{n+1})} \geq \sqrt{2+\varepsilon_1 r_n'(2)}$

$\Rightarrow r_n(0) \leq r_{n+1}(c_{n+1}) \leq r_{n+1}(2)$ or $r_n(0) \geq r_{n+1}(c_{n+1}) \geq r_{n+1}(2)$. Hence by induction result is true.

Hence the proof of lemma 1.1 is complete.

Proof of the Theorem 1:-

By Lemma 1.1 part (i), since the *sin* function is continuous and since the series $\sum_{i=1}^{n}\frac{\varepsilon_1\varepsilon_2...\varepsilon_i}{2^i}$

absolutely convergent, $r_n(2)$ and $r_n(0)$ is convergent. It is clear that

$$\lim_{n\to\infty} r_n(2) = \lim_{n\to\infty} r_n(0) = 2\sin\left(1+\sum_{i=1}^{\infty}\frac{\varepsilon_1\varepsilon_2...\varepsilon_i}{2^i}\right)\frac{\pi}{4}.$$

By lemma 1.1 part (ii) $\lim_{n\to\infty} r_n(c_n) = 2\sin\left(1 + \sum_{i=1}^{\infty} \frac{\varepsilon_1\varepsilon_2...\varepsilon_i}{2^i}\right)\frac{\pi}{4}$. Complete the proof of Theorem 1.

Notation:- We use the notation $\sqrt{2+\varepsilon_1\sqrt{2+\varepsilon_2\sqrt{2+...}}}$ to denote the $\lim_{n\to\infty} r_n(0) = \lim_{n\to\infty} r_n$.

**Theorem 2**

Let $(\varepsilon_n) \in \{-1,1\}^{\mathbb{N}}$ and $(\varepsilon_n)$ has the following property $\exists\, k \in \mathbb{N}$ s.t $\varepsilon_k = \varepsilon_{k+1} = -1$ & $\varepsilon_i = 1\ \forall i \geq k+2$, let $A \subset \{-1,1\}^{\mathbb{N}}$ be the set of all such sequences.

Define $\psi : \{-1,1\}^{\mathbb{N}} \setminus A \to [0,2]$ such that $\psi((\varepsilon_n)) = \lim_{n\to\infty} r_n$, where $(r_n)$ is the corresponding continued radical sequence. Then $\psi$ is *1-1* and on to.

Proof of Theorem 2

Observe that by Theorem 1, $\psi$ is well defined. To show $\psi$ is *1-1* suppose there exist two different sequences $(\varepsilon_n)$ and $(\varepsilon_n')$ in $\{-1,1\}^{\mathbb{N}} \setminus A$ such that $\psi((\varepsilon_n)) = \psi((\varepsilon_n'))$

$\Rightarrow \sqrt{2+\varepsilon_1\sqrt{2+\varepsilon_2\sqrt{2+...}}} = \sqrt{2+\varepsilon_1'\sqrt{2+\varepsilon_2'\sqrt{2+...}}}$. Let $i = \min\{k\,|\,\varepsilon_k \neq \varepsilon_k'\}$.

Then $\varepsilon_i\sqrt{2+\varepsilon_{i+1}\sqrt{2+...}} = \varepsilon_i'\sqrt{2+\varepsilon_{i+1}'\sqrt{2+...}}$, but since $\varepsilon_i \neq \varepsilon_i'$

$\sqrt{2+\varepsilon_{i+1}\sqrt{2+...}} = \sqrt{2+\varepsilon_{i+1}'\sqrt{2+...}} = 0$. But observe that for any $(\varepsilon_n) \in \{-1,1\}^{\mathbb{N}}$

$\lim_{n\to\infty} r_n = 0 \Leftrightarrow \varepsilon_1 = -1$ & $\varepsilon_i = 1\ \forall i \geq 2 \Rightarrow \varepsilon_{i+1} = \varepsilon_{i+1}' = -1$ & $\varepsilon_k = \varepsilon_k' = 1\ \forall k \geq i+2$.

But $\varepsilon_i \neq \varepsilon_i'\ \therefore \psi$ is $1-1 \Rightarrow (\varepsilon_n) \in A$ or $(\varepsilon_n') \in A$, contradiction $\therefore \psi$ is $1-1$.

To show $\psi$ is on to first we prove the following lemma.

**Lemma 2.1**

If $0 < \phi < \frac{\pi}{2}$ and $\cos 2^n\phi \neq 0\ \forall n \in \mathbb{N}$ then for any $n \in \mathbb{N}$,

$2\cos\phi = \sqrt{2+\varepsilon_1\sqrt{2+...+\varepsilon_n\sqrt{2+2\cos 2^{n+1}\phi}}} = r_n(2\cos 2^{n+1}\phi)$, where $\varepsilon_n = \frac{\cos 2^n\phi}{|\cos 2^n\phi|}$.

Proof of lemma 2.1

$2\cos\phi = \sqrt{2+2\cos 2\phi} = \sqrt{2+\dfrac{\cos 2\phi}{|\cos 2\phi|}\sqrt{2+2\cos 2^2\phi}}$. Observe that sign of the square root decide by sign of $\cos 2\phi$, it is equal to $\dfrac{\cos 2\phi}{|\cos 2\phi|}$. Hence, when $n=1$ the result is true.

Suppose for some $n \in \mathbb{N}$ the result is true. Observe that $2\cos 2^{n+1}\phi = \dfrac{\cos 2^{n+1}\phi}{|\cos 2^{n+1}\phi|}\sqrt{2+2\cos 2^{n+2}\phi}$.

Hence by induction result is true for any $n \in \mathbb{N}$. End of the proof of lemma 2.1.

To prove $\psi$ is on to first observe that when $\varepsilon_n = 1\ \forall n \in \mathbb{N}$, then $\psi((\varepsilon_n)) = 2$ and when $\varepsilon_1 = -1\ \&\ \varepsilon_n = 1\ \forall n \geq 2$, then $\psi((\varepsilon_n)) = 0$. In both cases $(\varepsilon_n) \notin A$.

Now let $x \in (0,2)$. Then there exist $\phi$ such that $0 < \phi < \dfrac{\pi}{2}$ and $x = 2\cos\phi$.

Suppose $\cos 2^n\phi \neq 0\ \forall n \in \mathbb{N}$. Then by above lemma $x = r_n(2\cos 2^{n+1}\phi)$ for any $n \in \mathbb{N}$.

where $\varepsilon_n = \dfrac{\cos 2^n\phi}{|\cos 2^n\phi|}$. By Theorem 1, $x = \lim\limits_{n\to\infty} r_n(2\cos 2^{n+1}\phi) = \lim\limits_{n\to\infty} r_n(0) = \sqrt{2+\varepsilon_1\sqrt{2+\varepsilon_2\sqrt{2+\ldots}}}$

$\therefore \psi\left(\left(\dfrac{\cos 2^n\phi}{|\cos 2^n\phi|}\right)\right) = x$. Also observe that $\left(\dfrac{\cos 2^n\phi}{|\cos 2^n\phi|}\right) \notin A$.

Now let $x \in (0,2)$ and $\exists k \in \mathbb{N}$ such that $\cos 2^k\phi = 0$. Let $i = \min\{k\,|\cos 2^k\phi = 0\}$.

Then $x = \sqrt{2+\varepsilon_1\sqrt{2+\ldots+\varepsilon_{i-1}\sqrt{2}}}$, where $\varepsilon_j = \dfrac{\cos 2^j\phi}{|\cos 2^j\phi|}$ for $j = 1,2,\ldots k-1$.

Define $(\varepsilon_n)$ as $\varepsilon_n = \dfrac{\cos 2^n\phi}{|\cos 2^n\phi|}$ when $1 \leq n \leq i-1$ and $\varepsilon_k = 1$, $\varepsilon_{k+1} = -1\ \&\ \varepsilon_n = 1\ \forall n \geq k+2$.

Then it is clear that $(\varepsilon_n) \notin A$ and $\psi((\varepsilon_n)) = x$. End of the proof of Theorem 2.

**Corollary 2.2**

Every $r \in [0,2]$ can be represented in the form of continued radical $r = \sqrt{2+\varepsilon_1\sqrt{2+\varepsilon_2\sqrt{2+...}}}$.

This represent is unique unless if we allowed using sequences $(\varepsilon_n)$ s.t. tail of these sequences consists of two consecutive -1's follows with infinitely many 1's.

Proof:- Clear.

**Theorem 3**

Let $x \in (0,2)$ & $k \in \mathbb{N}$, then $\exists\ \varepsilon_1, \varepsilon_2,...\varepsilon_k \in \{1,-1\}$ s.t. $x = \sqrt{2+\varepsilon_1\sqrt{2+...+\varepsilon_k\sqrt{2}}}$

$\Leftrightarrow x = 2\cos\left(\dfrac{\beta\pi}{2^{k+2}}\right)$, where $\beta$ is odd and $1 \leq \beta \leq 2^{k+1}-1$.

First we prove the following lemma.

**Lemma 3.1**

$\alpha$ is an odd integer and $-2^k + 1 \leq \alpha \leq 2^k - 1 \Leftrightarrow \exists\ \varepsilon_1, \varepsilon_2,...\varepsilon_k \in \{1,-1\}$ s.t. $\alpha = \sum_{i=1}^{k} 2^{i-1}\varepsilon_i$

Proof of lemma 3.1

Observe that for any $\varepsilon_1, \varepsilon_2,...\varepsilon_k \in \{1,-1\}$, $\sum_{i=1}^{k} 2^{i-1}\varepsilon_i$ is odd and $-2^k + 1 \leq \sum_{i=1}^{k} 2^{i-1}\varepsilon_i \leq 2^k - 1$.

Let $\sum_{i=1}^{k}(\varepsilon_i - \varepsilon_i')2^{i-1} = 0$. Suppose $\exists j \in \{1,2,...k\}$ s.t. $\varepsilon_i \neq \varepsilon_i'$. Let $r = \max\{j\ |\ \varepsilon_j \neq \varepsilon_j'\}$.

$\therefore \sum_{i=1}^{r}(\varepsilon_i - \varepsilon_i')2^{i-1} = 0$ and $\left|(\varepsilon_r - \varepsilon_r')\right| = 2$

$\Rightarrow (\varepsilon_r - \varepsilon_r')2^{r-1} = \sum_{i=1}^{r-1}(\varepsilon_i' - \varepsilon_i)2^{i-1} = 0 \Rightarrow \left|(\varepsilon_r - \varepsilon_r')\right|2^{r-1} \leq \sum_{i=1}^{r-1} 2 \cdot 2^{i-1} \Rightarrow 2^r \leq 2^r - 2$, conttradiction, $\therefore$

$\sum_{i=1}^{k}(\varepsilon_i - \varepsilon_i')2^{i-1} = 0 \Rightarrow \varepsilon_i = \varepsilon_i'$ for $i = 1,2,...k$, . Hence the representation $\sum_{i=1}^{k} 2^{i-1}\varepsilon_i$ is unique. $\therefore$ the set $\{(\varepsilon_1, \varepsilon_2,...\varepsilon_k)\ |\ \varepsilon_i \in \{1,-1\}$ for $i = 1,2,...k\}$ has $2^k$ elements.

$\therefore \{\alpha\ |\ -2^k + 1 \leq \alpha \leq 2^k$ and $\alpha$ is odd$\} = \left\{\sum_{i=1}^{k} 2^{i-1}\varepsilon_i\ |\ \varepsilon_i \in \{1,-1\}$ for $i = 1,2,...k\right\}$.Hence the lemma.

Proof of theorem 3:-

Let $x = \sqrt{2+\varepsilon_1\sqrt{2+...+\varepsilon_k\sqrt{2}}}$. Then $\varepsilon_1, \varepsilon_2,...\varepsilon_k$ can be extended to the sequence $(\varepsilon_n) \in \{1,-1\} \setminus A$ s.t. first k terms of the sequence is $\varepsilon_1, \varepsilon_2,...\varepsilon_k$ and $\varepsilon_{k+1} = 1$, $\varepsilon_{k+2} = -1$, $\varepsilon_{k+j} = 1 \forall j \geq 3$. Then observe that $x = \sqrt{2+\varepsilon_1\sqrt{2+...}}$. Therefore by Theorem 1,

$$x = 2\cos\left(1 - \sum_{i=1}^{k}\frac{\varepsilon_1\varepsilon_2...\varepsilon_i}{2^i} + \frac{\delta_k}{2^{k+1}} - \sum_{i=2}^{\infty}\frac{\delta_k}{2^{k+i}}\right)\frac{\pi}{4}, \text{ where } \delta_k = \varepsilon_1\varepsilon_2...\varepsilon_k.$$

Observe that $\sum_{i=2}^{\infty}\frac{\delta_k}{2^{k+i}} = \frac{\delta_k}{2^{k+1}} \therefore x = 2\cos\left(1 - \sum_{i=1}^{k}\frac{\varepsilon_1\varepsilon_2...\varepsilon_i}{2^i}\right)\frac{\pi}{4} = 2\cos\left(1 - \frac{1}{2}\sum_{i=1}^{k}\frac{\delta_i}{2^{i-1}}\right)\frac{\pi}{4}$, where

$\delta_i = \varepsilon_1\varepsilon_2...\varepsilon_i$ for $i = 1,2,...k$. $\Rightarrow x = 2\cos\left(1 - \frac{1}{2^k}\left(\sum_{i=1}^{k}\delta_i 2^{k-i}\right)\right)\frac{\pi}{4}$. Let $\alpha = \sum_{i=1}^{k}\delta_i 2^{k-i}$. Then by

Lemma 3.1 $-2^k + 1 \leq \alpha \leq 2^k - 1$ and $\alpha$ is odd. $\therefore x = 2\cos\left(1 - \frac{\alpha}{2^k}\right)\frac{\pi}{4} = 2\cos\left(\frac{\beta\pi}{2^{K+2}}\right)$

where $1 \leq \beta \leq 2^{k+1} - 1$ and $\beta$ is odd. Now suppose that $\beta$ is odd and $x = 2\cos\left(\frac{\beta\pi}{2^{k+2}}\right)$. Then observe that by setting $\frac{\beta}{2^k} = 1 - \frac{\alpha}{2^k}$, $\exists \alpha$ s.t. $-2^k + 1 \leq \alpha \leq 2^k - 1$, $\alpha$ is odd and $x = 2\cos\left(1 - \frac{\alpha}{2^k}\right)\frac{\pi}{4}$. Then by lemma 3.1, $\alpha = \sum_{i=1}^{k} 2^{i-1}\delta_i$. But then we can select $\varepsilon_1, \varepsilon_2,...\varepsilon_k \in \{1,-1\}$ s.t. $\varepsilon_1 = \delta_1$, $\varepsilon_2 = \frac{\delta_2}{\varepsilon_1}$ and inductively when $\varepsilon_1, \varepsilon_2...\varepsilon_{r-1}$ already selected, select $\varepsilon_r = \frac{\delta_r}{\varepsilon_1\varepsilon_2...\varepsilon_{r-1}}$.

$\therefore x = 2\cos\left(1 - \sum_{i=1}^{k}\frac{\varepsilon_1\varepsilon_2...\varepsilon_i}{2^i}\right)\frac{\pi}{4} = \sqrt{2+\varepsilon_1\sqrt{2+...+\varepsilon_k\sqrt{2}}}$. Hence the theorem.

Note:- Observe that any finite nested radical representation of the form $\sqrt{2+\varepsilon_1\sqrt{2+...+\varepsilon_n\sqrt{2}}}$

Can be extended to continued radical representation in two ways, by putting either $\varepsilon_{n+1} = 1, \varepsilon_{n+2} = -1$,

$\varepsilon_{n+3} = \varepsilon_{n+4} = ... = 1$ or $\varepsilon_{n+1} = -1, \varepsilon_{n+2} = -1$, $\varepsilon_{n+3} = \varepsilon_{n+4} = ... = 1$.

### Definition

Consider the continued radical $\sqrt{2+\varepsilon_1\sqrt{2+...}}$. Then this continued radical(or the sequence $(\varepsilon_n)$) is said to be eventually periodic iff There exist $k \in \mathbb{N}$ and $p \in \mathbb{N}$ s.t. $\varepsilon_{k+i+mp} = \varepsilon_{k+i} \ \forall \ m \in \mathbb{N} \cup \{0\}$ and $\forall i \in \{0, 1,... p-1\}$. The above continued radical is said to be "totally periodic" (or periodic) Iff

$\varepsilon_{1+i+mp} = \varepsilon_{1+i}$ $\forall\, m \in \mathbb{N} \cup \{0\}$ and $\forall\, i \in \{0, 1, \ldots p-1\}$. In both cases $p$ is called a period of the continued radical (or the sequence $(\varepsilon_n)$).

$p$ is called the minimum period of the continued radical if $p'$ is a period of continued radical then $p \le p'$.

**Theorem 4**

Consider the continued radical $r = \sqrt{2 + \varepsilon_1 \sqrt{2 + \ldots}}$ then $r$ is totally periodic $\Leftrightarrow r = 2\cos q\pi$ for some $q \in \left(0, \dfrac{1}{2}\right) \cap \mathbb{Q}$ and $q$ has an odd denominator. (We remove the trivial case $\varepsilon_i = 1,\ \forall\, n \in \mathbb{N}$.)

Furthermore, minimal period of $r = \min\{d\,|\,2^d \equiv \pm 1 \bmod s,\ s$ is the denominator of $q\}$

Proof:-

Suppose $(\varepsilon_n)$ is totally periodic. Let $p$ be a period of $(\varepsilon_n)$. Observe that if $\varepsilon_1\varepsilon_2\ldots\varepsilon_p = 1$ then $p$ is a period of $(\delta_n)$, (Where $\delta_n = \varepsilon_1\varepsilon_2\ldots\varepsilon_n$). If $\varepsilon_1\varepsilon_2\ldots\varepsilon_p = -1$, then $p$ is not a period of $(\delta_n)$ but $2p$ is a period. Hence without loss of generality let $p$ be a period of $(\delta_n)$ and $\delta_p = 1$.

By using summation of geometric series, $\displaystyle\sum_{i=1}^{\infty} \dfrac{\delta_i}{2^i} = \left(\dfrac{\delta_1}{2} + \dfrac{\delta_2}{2^2} + \ldots + \dfrac{\delta_{p-1}}{2^{p-1}} + \dfrac{1}{2^p}\right)\dfrac{2^p}{2^p - 1}$

$\Rightarrow 1 - \displaystyle\sum_{i=1}^{\infty} \dfrac{\delta_i}{2^i} = \left(2^p - 1 - (2^{p-1}\delta_1 + 2^{p-2}\delta_2 + \ldots + 2\delta_{p-1} + 1)\right)\dfrac{1}{2^p - 1}$

Observe that $2^p - 1 - (2^{p-1}\delta_1 + 2^{p-2}\delta_2 + \ldots + 2\delta_{p-1} + 1) \equiv 0 \bmod 4$. Hence the denominator of

$q = \left(1 - \displaystyle\sum_{i=1}^{\infty} \dfrac{\delta_i}{2^i}\right)\dfrac{1}{4}$ is odd. $q$ gets its maximum when $\delta_i = -1,\ \forall\, i \in \mathbb{N}$ and gets its minimum when

$\delta_i = 1,\ \forall\, i \in \mathbb{N}$, hence $0 < q < \tfrac{1}{2}$. Now suppose $r = 2\cos q\pi$, $q \in (0, \tfrac{1}{2})$, $q = \tfrac{t}{s}$, $s$ is odd & $(t, s) = 1$.

Hence $2\cos\tfrac{t\pi}{s} = \sqrt{2 + \varepsilon_1\sqrt{2 + \ldots}}$ ., where $\varepsilon_i = \dfrac{\cos 2^i \tfrac{t\pi}{s}}{\left|\cos 2^i \tfrac{t\pi}{s}\right|}$, (since $s$ is odd, $\cos 2^i \tfrac{t\pi}{s} \ne 0$).

Since $(2, t) = 1$, $2^{\varphi(s)} \equiv 1 \bmod s \Rightarrow 2^{\varphi(s)} t \equiv t \bmod s \Rightarrow \cos 2^{1+j+m\varphi(s)} \tfrac{t\pi}{s} = \cos 2^{1+j} \tfrac{t\pi}{s}$.

$\Rightarrow \varepsilon_{1+j+m\varphi(s)} = \varepsilon_{1+j}$ $\forall\, m \in \mathbb{N} \cup \{0\}$ and $\forall\, j \in \{0, 1, \ldots \varphi(s) - 1\} \Rightarrow (\varepsilon_n)$ is totally periodic and $\varphi(s)$ is a period. (Need not be a minimal period).

Let $p = \min\{d \mid 2^d \equiv \pm 1 \bmod s\}$ then $p \neq \emptyset$, since $2^{\varphi(s)} \equiv 1 \bmod s$. Then $2^d \equiv \pm 1 \bmod s$

$\Rightarrow 2^{md} t \equiv \pm t \bmod s \Rightarrow \cos\left(2^{1+j+md} \frac{t\pi}{s}\right) = \cos\left(2^{1+j} \frac{t\pi}{s}\right)$

$\Rightarrow \cos\left(2^{1+j+md} \frac{t\pi}{s}\right) = \cos\left(2^{1+j} \frac{t\pi}{s}\right) \Rightarrow \varepsilon_{1+j+dm} = \varepsilon_{1+j} \,\forall m \in \mathbb{N} \cup \{0\}$ and $\forall j \in \{0,1,...,(d-1)\} \Rightarrow d$ is a period.

Now suppose $d_0 \leq d$ and $d_0$ is a period. $\Rightarrow \varepsilon_{1+j+md_0} = \varepsilon_{1+j} \,\forall m \in \mathbb{N} \cup \{0\}$ & $\forall j \in \{0,1,...,(d_0-1)\}$.

$\Rightarrow \cos\left(2^{1+j+md_0} \frac{t\pi}{s}\right) = \cos\left(2^{1+j} \frac{t\pi}{s}\right) \Rightarrow 2^{1+j+md_0} \frac{t\pi}{s} = 2k\pi \pm 2^{1+j} \frac{t\pi}{s}, k \in \mathbb{Z}. \Rightarrow 2^{1+j}(2^{md_0} \frac{t}{s} \pm \frac{t}{s}) = 2k$

$(2^{md_0} \pm 1)\frac{t}{s} \in \mathbb{Z} \Rightarrow 2^{md_0} \equiv \pm 1 \bmod s$. By minimality of $d$ implies $d_0 = d$.

Hence the minimal period of $r = \min\{d \mid 2^d \equiv \pm 1 \bmod s\}$.

**Corollary 4.1**

$r = \sqrt{2 + \varepsilon_1 \sqrt{2 + ...}}$ eventually periodic $\Leftrightarrow r = 2\cos q\pi$ for some $q \in \left(0, \frac{1}{2}\right) \cap \mathbb{Q}$.

Proof:- Suppose $r$ is eventually periodic, then $r = 2\cos\left(1 - \sum_{i=1}^{\infty} \frac{\delta_i}{2^i}\right)\frac{\pi}{4}$ and $\exists k \in \mathbb{N}$ s.t. $\delta_{k+r+mp} = \delta_{k+r}$ for $r \in \{0,1,...,(p-1)\}$ for $m \in \mathbb{N}$. Using summation of geometric series,

$\sum_{i=k}^{\infty} \frac{\delta_i}{2^i} = \frac{1}{2^k}\left(\sum_{i=0}^{p-1} \frac{\delta_{k+i}}{2^i}\right)\left(\frac{2^p}{2^p - 1}\right) \in \mathbb{Q}$. Hence $q = \left(1 - \sum_{i=1}^{\infty} \frac{\delta_i}{2^i}\right)\left(\frac{1}{4}\right) \in \mathbb{Q}$.

Now consider the infinite nested radical representation of $r = 2\cos q\pi, q \in \left(0, \frac{1}{2}\right) \cap \mathbb{Q}$. If the denominator of $q$ has no odd factor other than 1, then $q = \frac{1}{2^2}$ or $q = \frac{\beta}{2^{k+2}}, k \in \mathbb{N}, 1 \leq \beta \leq 2^{k+1} - 1$ and

$\beta$ is odd. Hence note to the theorem 3 continued radical representation of $r$ is eventually periodic. Now suppose the denominator of $q$ is of the form $2^k s, k \geq 1$ and $s$ is an odd number greater than 3. Let $t$ be the numerator of $q$, then $t$ is odd. Then by lemma 2.1,

$r = \sqrt{2 + \varepsilon_1 \sqrt{2 + ... + \varepsilon_{k-1}\sqrt{2 + \varepsilon_k . 2\cos\frac{t_0\pi}{s}}}}$ , $\varepsilon_i = \frac{\cos\frac{t\pi}{2^{k-i}s}}{\left|\cos\frac{t\pi}{2^{k-i}s}\right|}$ , where $i \in \{1, 2, ..., (k-1)\}$.

$\varepsilon_k$ can be chosen as to set $0 < t_0 < \frac{s}{2}$. By theorem 4, $2\cos\frac{t_0\pi}{s}$ has totally periodic continued radical representation. Hence $r$ is eventually periodic.

Following definition is facilitating to develop a simple and practical method to find continued radical representation of $2\cos q\pi$ where $q$ is rational.

**Definition**

If $s$ is an odd number then semi order of 2 modulo $s = \min\{d \mid 2^d \equiv \pm 1 \bmod s\}$.

We give three steps to find continued radical representation of $2\cos q\pi$, where $q$ is rational.

Step1:- Write down $q$ as, $q = \frac{t}{2^k s}$, $k \in \mathbb{N} \cup \{0\}$, and $s$ and $t$ are relatively prime numbers. Without loss of generality we can assume that $0 < t < 2^{k-1} s$.

Step 2:- If $s = 1$, then $2\cos q\pi$ has the form $2\cos q\pi = \sqrt{2}$ or $2\cos q\pi = \sqrt{2 + \varepsilon_1 \sqrt{2 + \ldots \varepsilon_k \sqrt{2}}}$, where $\varepsilon_i = \dfrac{\cos 2^i q\pi}{|\cos 2^i q\pi|}$ to find $\varepsilon_i$ following practical method can be used.

Check whether the angle $2^i q\pi$ is in first, second, third or fourth quadrant. If it is in first or fourth quadrant $\varepsilon_i = 1$, if it is in second or third quadrant $\varepsilon_i = -1$.

Step 3:- If $s \geq 3$, then find the semi order of 2 modulo $s$, let it $p$. This is the minimal period of relevant continued radical. To find first $k$ coefficients $\varepsilon_1, \varepsilon_2, \ldots \varepsilon_k$ use the method in step 2. Now

$2\cos q\pi = \sqrt{2 + \varepsilon_1 \sqrt{2 + \ldots \varepsilon_k \sqrt{2 + \varepsilon_{k+1}.2\cos \frac{t_0 \pi}{3}}}}$, where $\varepsilon_{k+1}$ can be selected as to set $0 < t_0 < \frac{s}{2}$. Again, use the method in step 2 to find first $p$ coefficients of relevant continued radical of $2\cos \frac{t_0 \pi}{s}$. Now this is totally periodic with periodic block $(\varepsilon_{k+2}, \varepsilon_{k+3}, \ldots, \varepsilon_{k+p+1})$.

Example:- Find the continued radicals representations of $2\cos \frac{21\pi}{136}$.

Observe that the semi order of 2 modulo 17 is 4, since $\min\{d \mid 2^d \equiv \pm 1 \bmod 17\} = 4$.

$2\cos \frac{2.21\pi}{2^3.17} \in 1^{st}$ quadrant, $2\cos \frac{2^2.21\pi}{2^3.17} \in 2^{nd}$ quadrant, $2\cos \frac{2^3.21\pi}{2^3.17} \in 3^{rd}$ quadrant.

$\Rightarrow \varepsilon_1 = 1, \varepsilon_2 = -1, \varepsilon_3 = -1$.

$\therefore 2\cos \frac{21\pi}{136} = \sqrt{2 + \sqrt{2 - \sqrt{2 - \sqrt{2 + 2\cos \frac{21\pi}{17}}}}} = \sqrt{2 + \sqrt{2 - \sqrt{2 - \sqrt{2 + \varepsilon_4.2\cos \frac{4\pi}{17}}}}}$, where $\varepsilon_4 = -1$.

Now calculate first four coefficients of $2\cos \frac{4\pi}{17}$.

$2\cos\frac{2.4\pi}{17} \in 1^{st}$ quadrant, $2\cos\frac{2^2.4\pi}{17} \in 2^{nd}$ quadrant, $2\cos\frac{2.4\pi}{17} \in 1^{st}$ quadrant, $2\cos\frac{2^3.4\pi}{17} \in 4^{th}$ quadrant,

$2\cos\frac{2^4.4\pi}{17} \in 4^{th}$ quadrant, $\Rightarrow \varepsilon_5 = 1, \varepsilon_6 = -1, \varepsilon_7 = 1, \varepsilon_8 = 1$. Hence $2\cos\frac{4\pi}{17}$ has repeated block $(+,-,+,+)$. Consequently representation of $2\cos\frac{21\pi}{136}$ is $+--,-,+-++,+-++,+-++,....$

## Section 2

**Theorem 5**

Let $(r_n)$ be a totally periodic continued radical sequence of minimal period $p$, with the corresponding block $(\varepsilon_1, \varepsilon_2, ...\varepsilon_p)$. Let $n = mp + j$, where $0 \leq j \leq p-1$, $m \in \mathbb{N} \cup \{0\}$. Then,

$$r_n = r_{mp+j} = \begin{cases} 2\cos\left(1 - \frac{2^p \sigma_p}{2^p-1} + \frac{1}{2^{mp}} \frac{2^p \sigma_p}{(2^p-1)} - \frac{\sigma_j}{2^{mp}}\right)\frac{\pi}{4}, & \text{when } \delta_p = 1 \\ 2\cos\left(1 - \frac{2^p \sigma_p}{2^p+1} + \frac{(-1)^m}{2^{mp}} \frac{2^p \sigma_p}{(2^p+1)} - \frac{(-1)^m \sigma_j}{2^{mp}}\right)\frac{\pi}{4}, & \text{when } \delta_p = -1 \end{cases}$$

Where $\sigma_j = \sum\limits_{i=0}^{j} \frac{\delta_i}{2^i}$, $\delta_0 = 0$, $\delta_i = \varepsilon_1 \varepsilon_2 ... \varepsilon_i$.

Proof:-

By lemma 1.1, $r_n = 2\cos\left(1 - \sum\limits_{i=1}^{n} \frac{\varepsilon_1 \varepsilon_2 ... \varepsilon_i}{2^i}\right)\frac{\pi}{4} = 2\cos\left(1 - \sum\limits_{i=1}^{n} \frac{\delta_i}{2^i}\right)\frac{\pi}{4}$.

Observe that $\sum\limits_{i=1}^{n} \frac{\delta_i}{2^i} = \sum\limits_{i=1}^{mp} \frac{\delta_i}{2^i} + \frac{\delta_{mp}}{2^{mp}} \sum\limits_{i=1}^{j} \frac{\delta_i}{2^i}$, but $\delta_{mp} = \underbrace{(\varepsilon_1\varepsilon_2...\varepsilon_p)(\varepsilon_1\varepsilon_2...\varepsilon_p)...(\varepsilon_1\varepsilon_2...\varepsilon_p)}_{m \text{ blocks}} = (\delta)^{pm}$.

$\therefore \sum\limits_{i=1}^{n} \frac{\delta_i}{2^i} = \sum\limits_{i=1}^{mp} \frac{\delta_i}{2^i} + \frac{(\delta_p)^m}{2^{mp}} \sum\limits_{i=1}^{j} \frac{\delta_i}{2^i} \Rightarrow \sigma_n = \sigma_{mp} + \frac{(\delta_p)^m}{2^{mp}} \sigma_j$. Also observe that,

$$\sigma_{mp} = \left(\frac{\delta_1}{2} + \frac{\delta_2}{2^2} + ... + \frac{\delta_p}{2^p}\right) + \delta_p\left(\frac{\delta_1}{2^{p+1}} + \frac{\delta_2}{2^2} + ... + \frac{\delta_p}{2^{2p}}\right) + ... + (\delta_p)^{m-1}\left(\frac{\delta_1}{2^{(m-1)p+1}} + \frac{\delta_2}{2^{(m-1)p+2}} + ... + \frac{\delta_p}{2^{mp}}\right).$$

$$\sigma_{mp} = \begin{cases} \frac{2^p \sigma_p}{2^p-1} + \frac{1}{2^{mp}} \frac{2^p \sigma_p}{(2^p-1)}, & \text{when } \delta_p = 1 \\ \frac{2^p \sigma_p}{2^p+1} - \frac{(-1)^m}{2^{mp}} \frac{2^p \sigma_p}{(2^p+1)} & \text{when } \delta_p = -1 \end{cases}$$

Hence the theorem.

**Corollary 5.1**

$$\lim_{n\to\infty} r_n = \begin{cases} 2\cos\left(1 - \frac{2^p \sigma_p}{2^p - 1}\right)\frac{\pi}{4}, & \text{when } \delta_p = 1 \\ 2\cos\left(1 - \frac{2^p \sigma_p}{2^p - 1}\right)\frac{\pi}{4}, & \text{when } \delta_p = -1 \end{cases}$$

Proof :- Take the limit of both sides in the identity of theorem 5.

**Theorem 6**

Let $(r_n)$ be a totally periodic continued radical sequence of 2. Let $(u_n)$ be the sequence defined as,

$u_n = 2^n (r_\infty - r_n)$, where $r_\infty = \lim_{n\to\infty} r_n$.

Case I: -

If $\delta_p = 1$ then $(u_n)$ has $p$ distinct limits points which are, $2^j \left(\frac{2^p \sigma_p}{2^p - 1} - \sigma_j\right) \frac{\pi}{2} \sin\left(1 - \frac{2^p \sigma_p}{2^p - 1}\right)\frac{\pi}{4}$ for $j = 0, 1, 2 \ldots (p-1)$.

Case II:-

If $\delta_p = -1$ then $(u_n)$ has $2p$ distinct limits points which are, $2^j \left(\frac{2^p \sigma_p}{2^p + 1} - \sigma_j\right) \frac{\pi}{2} \sin\left(1 - \frac{2^p \sigma_p}{2^p + 1}\right)\frac{\pi}{4}$ for $j = 0, 1, 2 \ldots (2p-1)$.

Proof of the theorem:-

Let $\delta_p = 1$, the sequence $(r_n)$ can be partitioned in to $p$ subsequences according to $(r_{mp+j})$

$j = 0, 1, 2 \ldots (p-1)$, $r_o = \sqrt{2}$, $m \in \mathbb{N} \cup \{0\}$. Since the sequence $(r_n)$ converges, each of these subsequences convergent to $r_\infty = 2\cos\left(1 - \frac{2^p \sigma_p}{2^p - 1}\right)\frac{\pi}{4}$. Consequently the sequence $(u_n)$ can be partitioned in to $p$ subsequences as $(u_{mp+j})$ $j = 0, 1, 2 \ldots (p-1)$, $m \in \mathbb{N} \cup \{0\}$.

Consider the subsequence $(u_{mp+j})$ for some $j = 0, 1, 2 \ldots (p-1)$ & $m \in \mathbb{N} \cup \{0\}$.

From theorem 5, $u_{mp+j} = 2^{mp+j+1}\left\{\cos\left(1 - \frac{2^p \sigma_p}{2^p - 1}\right)\frac{\pi}{4} - \cos\left(1 - \frac{2^p \sigma_p}{2^p - 1} + \frac{1}{2^{mp}} \frac{2^p \sigma_p}{(2^p - 1)} - \frac{\sigma_j}{2^{mp}}\right)\frac{\pi}{4}\right\}$.

Now using the cosine addition formula $\cos A - \cos B = 2\sin\frac{(A+B)}{2}\sin\frac{(B-A)}{2}$,

$$u_{mp+j} = 2^{mp+j+2}.\sin\left(\frac{1}{2^{mp}}.\frac{2^p\sigma_p}{2^p-1} - \frac{\sigma_j}{2^{mp}}\right)\frac{\pi}{8}.\sin\left(1 - \frac{2^p\sigma_p}{2^p-1} + \frac{1}{2^{mp+1}}\frac{2^p\sigma_p}{(2^p-1)} - \frac{\sigma_j}{2^{mp+1}}\right)\frac{\pi}{4}$$

$$\therefore \lim_{m\to\infty} u_{mp+j} = \lim_{m\to\infty} 2^{mp+j+2}.\sin\left(\frac{2^p\sigma_p}{2^p-1} - \sigma_j\right)\frac{\pi}{8.2^{mp}}.\lim_{m\to\infty}\sin\left(1 - \frac{2^p\sigma_p}{2^p-1} + \frac{1}{2^{mp+1}}\frac{2^p\sigma_p}{(2^p-1)} - \frac{\sigma_j}{2^{mp+1}}\right)\frac{\pi}{4}$$

$$= \lim_{m\to\infty} 2^{j+2}.\frac{\left(\frac{2^p\sigma_p}{2^p-1} - \sigma_j\right)\frac{\pi}{8}.\sin\left(\frac{2^p\sigma_p}{2^p-1} - \sigma_j\right)\frac{\pi}{8.2^{mp}}}{\left(\frac{2^p\sigma_p}{2^p-1} - \sigma_j\right)\frac{\pi}{8.2^{mp}}}\sin\left(1 - \frac{2^p\sigma_p}{2^p-1}\right)\frac{\pi}{4}, \text{ \{using } \lim_{k\to\infty} k\sin\frac{\alpha}{k} = \alpha\text{\}}$$

$$= 2^j\left(\frac{2^p\sigma_p}{2^p-1} - \sigma_j\right)\frac{\pi}{2}\sin\left(1 - \frac{2^p\sigma_p}{2^p-1}\right)\frac{\pi}{4}.$$

Different subsequences give different limits for $j = 0,1,2...(p-1)$, since $\sigma_i \neq \sigma_j$ for $i \neq j$.

Particularly, the sub sequence $u_{mp}$ converges to $\left(\frac{2^p\sigma_p}{2^p-1} - \sigma_j\right)\frac{\pi}{2}\sin\left(1 - \frac{2^p\sigma_p}{2^p-1}\right)\frac{\pi}{4}$.

Now let $\delta_p = -1$, partition the sequence $(r_n)$ in to 2p different subsequences accordance to $(r_{2mp+j})$ for $j = 0,1,2...(2p-1)$, $r_o = \sqrt{2}$, $m \in \mathbb{N} \cup \{0\}$. Since the sequence $(r_n)$ converges, each of these subsequences convergent to $r_\infty = 2\cos\left(1 - \frac{2^p\sigma_p}{2^p+1}\right)\frac{\pi}{4}$. Consequently the sequence $(u_n)$ can be partitioned in to 2p subsequences as $(u_{2mp+j})$ $j = 0,1,2...(2p-1)$, $m \in \mathbb{N} \cup \{0\}$.

From theorem 5, substituting 2m instead of m, we can get the identity,

$$r_{2mp+j} = 2\cos\left(1 - \frac{2^p\sigma_p}{2^p+1} + \frac{1}{2^{2mp}}\frac{2^p\sigma_p}{(2^p+1)} - \frac{\sigma_j}{2^{2mp}}\right)\frac{\pi}{4}, \text{ for } j = 0,1,2...(2p-1)$$

Consider the subsequence $(u_{2mp+j})$ for some $j = 0,1,2...(2p-1)$ & $m \in \mathbb{N} \cup \{0\}$.

From theorem 5, $u_{2mp+j} = 2^{mp+j+1}\left\{\cos\left(1 - \frac{2^p\sigma_p}{2^p+1}\right)\frac{\pi}{4} - \cos\left(1 - \frac{2^p\sigma_p}{2^p+1} + \frac{1}{2^{2mp}}\frac{2^p\sigma_p}{(2^p+1)} - \frac{\sigma_j}{2^{2mp}}\right)\frac{\pi}{4}\right\}$

$$= 2^{2mp+j+2}.\sin\left(\frac{1}{2^{2mp}}.\frac{2^p\sigma_p}{2^p+1} - \frac{\sigma_j}{2^{2mp}}\right)\frac{\pi}{8}.\sin\left(1 - \frac{2^p\sigma_p}{2^p+1} + \frac{1}{2^{mp+1}}\frac{2^p\sigma_p}{(2^p+1)} - \frac{\sigma_j}{2^{mp+1}}\right)\frac{\pi}{4}.$$

By similar calculations as in case I,

$$\lim_{m\to\infty} u_{2mp+j} = 2^j\left(\frac{2^p\sigma_p}{2^p+1} - \sigma_j\right)\frac{\pi}{2}\sin\left(1 - \frac{2^p\sigma_p}{2^p+1}\right)\frac{\pi}{4}.$$

Observe that $\sigma_i \neq \sigma_j$ for $i \neq j$, $i, j \in \{0,1,2...(2p-1)\}$. Hence Different subsequences give different limits for $j = 0,1,2...(2p-1)$. Hence the theorem.

**Corollary 6.1**

Limits points of the sequence $u_n = 2^n (r_\infty - r_n)$, where $r_\infty = \lim\limits_{n \to \infty} r_n$, have the form $\alpha\pi$, where $\alpha$ is an algebraic number.

Proof:- If $q$ is rational then $\sin q\pi$ is algebraic[11]. Hence each limit point is of the form $\alpha\pi$, where $\alpha$ is an algebraic number.

**Result of Theorem 6**

Let $(r_n)$ be a totally periodic continued radical sequence of 2 with minimal period 3 with periodic block $(-,+,-) \Rightarrow \varepsilon_1 = -1, \varepsilon_2 = 1, \varepsilon_3 = -1 \Rightarrow \delta_1 = -1, \delta_2 = -1, \delta_3 = 1$ ($\& \delta_0 = 0$)

$$\Rightarrow \sigma_1 = \sum_{i=0}^{1} \tfrac{\delta_i}{2^i} = -\tfrac{1}{2}, \ \sigma_2 = \sum_{i=0}^{2} \tfrac{\delta_i}{2^i} = -\tfrac{1}{2} - \tfrac{1}{2^2} = -\tfrac{3}{4}, \ \sigma_3 = \sum_{i=0}^{3} \tfrac{\delta_i}{2^i} = -\tfrac{1}{2} - \tfrac{1}{2^2} + \tfrac{1}{2^3} = -\tfrac{5}{8}, (\& \sigma_0 = 0).$$

$\Rightarrow \lim\limits_{n \to \infty} r_n = r_\infty = 2\cos\left(1 - \tfrac{2^3 \sigma_3}{2^3 - 1}\right)\tfrac{\pi}{4} = 2\cos\tfrac{3\pi}{7}$. Hence by theorem 6 (case I) the sequence

$$2^n \left( 2\cos\tfrac{3\pi}{7} - \underbrace{\sqrt{2 - \sqrt{2 + \sqrt{2 - ... \sqrt{2}}}}}_{n \text{ roots}} \right) \text{ has three limit points, which are },$$

$2^0 \left( \tfrac{2^3 \sigma_3}{2^3 - 1} - \sigma_0 \right) \tfrac{\pi}{2} \sin \tfrac{3\pi}{7} = -\tfrac{5\pi}{14} \sin \tfrac{3\pi}{7}$, $2^1 \left( \tfrac{2^3 \sigma_3}{2^3 - 1} - \sigma_1 \right) \tfrac{\pi}{2} \sin \tfrac{3\pi}{7} = -\tfrac{3\pi}{14} \sin \tfrac{3\pi}{7}$,

$2^2 \left( \tfrac{2^3 \sigma_3}{2^3 - 1} - \sigma_1 \right) \tfrac{\pi}{2} \sin \tfrac{3\pi}{7} = \tfrac{\pi}{14} \sin \tfrac{3\pi}{7}$.

Particularly $\lim\limits_{n \to \infty} 2^{3n} \left( 2\cos\tfrac{3\pi}{7} - \underbrace{\sqrt{2 - \sqrt{2 + \sqrt{2 - ... \sqrt{2}}}}}_{3n \text{ roots}} \right) = -\tfrac{5\pi}{14} \sin \tfrac{3\pi}{7}$.

Let $(r_n)$ be a totally periodic continued radical sequence of 2 with minimal period 1 with repeated sign $(-)$ i.e. $r_n = \underbrace{\sqrt{2 - \sqrt{2 - \sqrt{2 - ... \sqrt{2}}}}}_{n \text{ roots}} \Rightarrow \varepsilon_1 = -1, \varepsilon_2 = -1 \Rightarrow \delta_1 = -1, \delta_2 = 1, (\& \delta_0 = 0)$

$\Rightarrow \sigma_1 = \sum_{i=0}^{1} \frac{\delta_i}{2^i} = -\frac{1}{2} \Rightarrow \lim_{n\to\infty} r_n = r_\infty = 2\cos\left(1 - \frac{2\sigma_1}{2+1}\right)\frac{\pi}{4} = 2\cos\frac{\pi}{3}$. Hence by theorem 6(case II) the sequence $2^n\left(1 - \sqrt{2 - \underbrace{\sqrt{2 - \sqrt{2 - ... - \sqrt{2}}}}_{n \text{ roots}}}\right)$ has two limit points, which are,

$2^0\left(\frac{2\sigma_1}{2+1} - \sigma_0\right)\frac{\pi}{2}\sin\frac{\pi}{3} = -\frac{\sqrt{3}\pi}{12}$ and $2^1\left(\frac{2\sigma_1}{2+1} - \sigma_1\right)\frac{\pi}{2}\sin\frac{\pi}{3} = \frac{\sqrt{3}\pi}{12}$. Particularly the sequence

$2^{2n}\left(1 - \sqrt{2 - \underbrace{\sqrt{2 - \sqrt{2 - ... - \sqrt{2}}}}_{2n \text{ roots}}}\right)$ converges to $-\frac{\sqrt{3}\pi}{12}$ and the sequence

$2^{2n-1}\left(1 - \sqrt{2 - \underbrace{\sqrt{2 - \sqrt{2 - ... - \sqrt{2}}}}_{2n-1 \text{ roots}}}\right)$ converges to $\frac{\sqrt{3}\pi}{12}$.

Let $(r_n)$ be a totally periodic continued radical sequence of 2 with periodic block $(+, -)$.

$\Rightarrow \varepsilon_1 = -1, \varepsilon_2 = -1 \Rightarrow \delta_1 = 1, \delta_2 = -1, \delta_3 = -1, \delta_4 = 1 \ (\& \delta_0 = 0)$

$\Rightarrow \sigma_1 = \sum_{i=0}^{1} \frac{\delta_i}{2^i} = \frac{1}{2}, \ \sigma_2 = \sum_{i=0}^{2} \frac{\delta_i}{2^i} = \frac{1}{2} - \frac{1}{2^2} = \frac{1}{4}, \sigma_3 = \sum_{i=0}^{3} \frac{\delta_i}{2^i} = \frac{1}{2} - \frac{1}{2^2} - \frac{1}{2^3} = \frac{1}{8} (\& \sigma_0 = 0)$.

$\Rightarrow \lim_{n\to\infty} r_n = r_\infty = 2\cos\left(1 - \frac{2^2 \sigma_3}{2^2+1}\right)\frac{\pi}{4} = 2\cos\frac{\pi}{5}$. Hence by theorem 6 (case II) the sequence

$2^n\left(2\cos\frac{\pi}{5} - r_n\right)$ has four limit points, which are $2^0\left(\frac{2^2\sigma_2}{2^2+1} - \sigma_0\right)\frac{\pi}{2}\sin\frac{\pi}{5} = \frac{\pi}{10}\sin\frac{\pi}{5}$,

$2^1\left(\frac{2^2\sigma_2}{2^2+1} - \sigma_1\right)\frac{\pi}{2}\sin\frac{\pi}{5} = -\frac{3\pi}{10}\sin\frac{\pi}{5}$, $2^2\left(\frac{2^2\sigma_2}{2^2+1} - \sigma_2\right)\frac{\pi}{2}\sin\frac{\pi}{5} = -\frac{\pi}{10}\sin\frac{\pi}{5}$, $2^2\left(\frac{2^2\sigma_2}{2^2+1} - \sigma_3\right)\frac{\pi}{2}\sin\frac{\pi}{5} = \frac{3\pi}{10}\sin\frac{\pi}{5}$.

Particularly $\lim_{n\to\infty} 2^{2n}\left(2\cos\frac{\pi}{5} - \sqrt{2 + \underbrace{\sqrt{2 - \sqrt{2 - ... - \sqrt{2}}}}_{2n \text{ roots}}}\right) = \frac{\pi}{10}\sin\frac{\pi}{5}$.

**Theorem 7**

Let $q \in \mathbb{Q}$, such that $0 < q < \frac{1}{2}$ with denominator of $q$ is odd, i.e. $q = \frac{t}{2k+1}$ for some $k \in \mathbb{N}$ and $t \in \{1, 2, 3...k\}$. Let $p$ be the minimal period of corresponding infinite nested radical of $2\cos q\pi$. i.e. $p = \min\{d \mid 2^d \equiv \pm 1 \mod (2k+1)\}$.

Define the sequence $(s_n)$ as follows.

$s_0 = \sqrt{2}$, $s_{pn+i} = \sqrt{2 + \varepsilon_{p+1-i} s_{pn+i-1}}$ for $n \in \mathbb{N} \cup \{0\}$ and $i \in \{1, 2, 3, \ldots p-1, p\}$,

where $\varepsilon_i = \dfrac{\cos 2^i q\pi}{|\cos 2^i q\pi|}$ for $i = 1, 2, 3, \ldots p$

Then, $\dfrac{1}{\pi} \dfrac{2\cos 2q\pi}{(1-4q)\sin q\pi} = \left(\dfrac{2\cos q\pi + \sqrt{2}}{2}\right) \prod_{i=0}^{\infty} \prod_{j=1}^{p} \left(\dfrac{|2\cos 2^{p-j} q\pi| + s_{pi+j}}{2}\right).$

Proof of theorem 7:-

Let $u_n = 2^n (2\cos q\pi - r_n)$, where $r_n$ is the $n^{th}$ iteration of the continued radical of $2\cos q\pi$, $q$ is as of the theorem. Then from theorem 6, noticed that the sequence $u_{pn}$ converges to $(1-4q)\frac{\pi}{2}\sin q\pi$.

(Simple substitution $q = \left(1 - \dfrac{2^p \sigma_p}{2^p - 1}\right)\frac{1}{4}$ or $q = \left(1 - \dfrac{2^p \sigma_p}{2^p + 1}\right)\frac{1}{4}$ gives this.)

By the definition of the sequence $(s_n)$ we can have the identity, for $1 \leq j \leq p$,

$|2\cos 2^{p-j} q\pi| - s_{pn+j} = \dfrac{(|2\cos 2^{p-j} q\pi| - s_{pn+j})(|2\cos 2^{p-j} q\pi| + s_{pn+j})}{(|2\cos 2^{p-j} q\pi| + s_{pn+j})}$

$= \dfrac{4\cos^2 2^{p-j} q\pi - 2 - \varepsilon_{p-j+1} s_{p\,nj-1}}{(|2\cos 2^{p-j} q\pi| + s_{pn+j})} = \dfrac{2\cos 2^{p-j+1} q\pi - \varepsilon_{p-j+1} s_{p\,nj-1}}{(|2\cos 2^{p-j} q\pi| + s_{pn+j})}$

Now since $\varepsilon_i = \dfrac{\cos 2^i q\pi}{|\cos 2^i q\pi|}$, observe that $2\cos 2^{p-j+1} q\pi$ and $\varepsilon_{p-j+1}$ have same sign. Therefore,

$|2\cos 2^{p-j} q\pi| - s_{pn+j} = \dfrac{\varepsilon_{p-j+1}(|2\cos 2^{p-j+1} q\pi| - s_{p\,nj-1})}{(|2\cos 2^{p-j} q\pi| + s_{pn+j})}$, for $1 \leq j \leq p$.

Now observe that by iterating above identity for $j = 1, 2, \ldots p$, we can derive the identity

$2\cos q\pi - s_{p(n+1)} = \dfrac{(\delta_p)^{n+1}(2\cos q\pi - s_0)}{\prod_{i=0}^{n}\prod_{j=1}^{p}(|2\cos 2^{p-j} q\pi| + s_{pi+j})}$. Now by the definition of sequences $(r_n)$ and $(s_n)$

notice that $r_{pn} = s_{pn}$ for all $n \in \mathbb{N}$. By theorem 6, $\lim_{n \to \infty} u_{p(n+1)} = (1-4q)\frac{\pi}{2}\sin q\pi$, where,

$$u_{p(n+1)} = 2^{p(n+1)}\left(2\cos q\pi - r_{p(n+1)}\right) = 2^{p(n+1)}\left(2\cos q\pi - s_{p(n+1)}\right) = \frac{2^{p(n+1)}\left(\delta_p\right)^{n+1}\left(2\cos q\pi - \sqrt{2}\right)}{\prod_{i=0}^{n}\prod_{j=1}^{p}\left(\left|2\cos 2^{p-j}q\pi\right| + s_{pi+j}\right)}$$

$$= \frac{\left(\delta_p\right)^{n+1}\left(2\cos q\pi - \sqrt{2}\right)}{\prod_{i=0}^{n}\prod_{j=1}^{p}\left(\frac{\left|2\cos 2^{p-j}q\pi\right| + s_{pi+j}}{2}\right)}$$

Now if $\delta_p = 1$, since the sequence $u_{p(n+1)}$ converges the product $\prod_{i=0}^{\infty}\prod_{j=1}^{p}\left(\frac{\left|2\cos 2^{p-j}q\pi\right| + s_{pi+j}}{2}\right)$ converges and taking limits in both sides, and with some manipulation gives the result.

If $\delta_p = -1$, substituting 2p, instead of $p$, we can have the following identity.

$$2^{2p(n+1)}\left(2\cos q\pi - s_{2p(n+1)}\right) = \frac{\left(2\cos q\pi - \sqrt{2}\right)}{\prod_{i=0}^{n}\prod_{j=1}^{2p}\left(\frac{\left|2\cos 2^{2p-j}q\pi\right| + s_{2pi+j}}{2}\right)}$$

But observe that, since $2^d \equiv -1 \mod(2k+1)$, $2^{2p-j} + 2^{p-j} \equiv 0 \mod(2k+1)$, where $2k+1$ is the denominator of $q$. Hence $\left|2\cos 2^{2p-j}q\pi\right| = \left|2\cos 2^{p-j}q\pi\right|$ for $j = 0, 1, 2, \ldots p$. This implies,

$$\prod_{i=0}^{n}\prod_{j=1}^{2p}\left(\frac{\left|2\cos 2^{2p-j}q\pi\right| + s_{2pi+j}}{2}\right) = \prod_{i=0}^{2n+1}\prod_{j=1}^{p}\left(\frac{\left|2\cos 2^{p-j}q\pi\right| + s_{pi+j}}{2}\right).$$ Now as in earlier, taking the limits in both sides and

with some manipulation gives $\frac{1}{\pi}\frac{2\cos 2q\pi}{(1-4q)\sin q\pi} = \left(\frac{2\cos q\pi + \sqrt{2}}{2}\right)\lim_{n\to\infty}\prod_{i=0}^{2n+1}\prod_{j=1}^{p}\left(\frac{\left|2\cos 2^{p-j}q\pi\right| + s_{pi+j}}{2}\right).$

Take $w_i = \prod_{j=1}^{p}\left(\frac{\left|2\cos 2^{p-j}q\pi\right| + s_{pi+j}}{2}\right)$, now since $\prod_{i=0}^{2n+1} w_i$ converges, $\lim_{n\to\infty} w_{2n}w_{2n+1} = 1$. But observe that the

definition of the sequence $(s_n)$ implies $\lim_{n\to\infty} s_{2np+j} = \lim_{n\to\infty} s_{2(n+1)+j}$ for $j = 0, 1, 2, \ldots p$.

$\therefore \lim_{n\to\infty} w_{2n} = \lim_{n\to\infty} w_{2n+1} = 1$. Hence $\prod_{i=0}^{n} w_i$ converges and $\lim_{n\to\infty}\prod_{i=0}^{n} w_i = \lim_{n\to\infty}\prod_{i=0}^{2n+1} w_i$. This gives the result,

$$\frac{1}{\pi}\frac{2\cos 2q\pi}{(1-4q)\sin q\pi} = \left(\frac{2\cos q\pi + \sqrt{2}}{2}\right)\prod_{i=0}^{\infty}\prod_{j=1}^{p}\left(\frac{\left|2\cos 2^{p-j}q\pi\right| + s_{pi+j}}{2}\right).$$ Hence the theorem.

## Results of theorem7

let $q = \dfrac{1}{3}$ . Then by theorem7,

$$\frac{2\sqrt{3}}{\pi} = \frac{1+\sqrt{2}}{2} \cdot \frac{1+\sqrt{2-\sqrt{2}}}{2} \cdot \frac{1+\sqrt{2-\sqrt{2-\sqrt{2}}}}{2} \cdots \quad \text{(Aaron Levin's result)}$$

Let $q = \dfrac{1}{5}$, Then by theorem7,

$$\frac{(\sqrt{5}-1)5\sqrt{2}}{\pi\sqrt{5-\sqrt{5}}} = \left(\frac{(\sqrt{5}+1)\frac{1}{2}+\sqrt{2}}{2}\right)\left(\frac{(\sqrt{5}-1)\frac{1}{2}+\sqrt{2-\sqrt{2}}}{2} \cdot \frac{(\sqrt{5}+1)\frac{1}{2}+\sqrt{2+\sqrt{2-\sqrt{2}}}}{2}\right)\left(\frac{(\sqrt{5}-1)\frac{1}{2}+\sqrt{2-\sqrt{2+\sqrt{2-\sqrt{2}}}}}{2} \cdot \frac{(\sqrt{5}+1)\frac{1}{2}+\sqrt{2+\sqrt{2-\sqrt{2+\sqrt{2-\sqrt{2}}}}}}{2}\right)\cdots$$